\documentclass[11pt,a4paper,geompsfi]{article}
\usepackage{amsfonts}

\usepackage[all]{xy}
\usepackage[latin1]{inputenc}        
\usepackage[dvips]{graphics}
\usepackage[dvips]{graphicx}
\usepackage{amssymb}
\usepackage{amsmath}
\usepackage{amstext}
\usepackage{amsbsy}
\usepackage{amsopn}
\usepackage{amsthm}
\usepackage{epsfig,rotating}
\usepackage{amscd}
\usepackage{amsxtra}
\usepackage{mathrsfs}
\usepackage{upref}

\let\small=\footnotesize
\def\longpage#1{\newdimen\addstr  \addstr=#1\baselineskip
  \advance\topmargin-2\baselineskip \advance\textheight4\baselineskip
  \advance\topmargin-0.5\addstr     \advance\textheight\addstr
  \advance\oddsidemargin-0.5\addstr \advance\textwidth\addstr}
\longpage{0.44}

       \def\R{\mathbb R}        
                    
                  \def\f{\varphi}       
          \def\:{\colon\,}        
               
\def\s{\sigma}                 \def\g{\gamma} 
\def\e{\varepsilon}            \def\G{\Gamma} 
                 
\def\sph{\mathbb S}

\def\k{\kappa}                 
\def\fin{\hfill{$\square$}}    
\def\d{\delta}

\def\N{\mathbb N}

\author{Ricardo \sc Uribe-Vargas \\ {\small Coll{\`e}ge de France, 11 Pl. Marcelin-Berthelot, 
75005 Paris.}\\ 
{\small uribe@math.jussieu.fr \ \ \ 
www.math.jussieu.fr/{\small $\sim$}uribe/ }}

\date\empty                     
\title{On Vertices,  Focal Curvatures and
Differential Geometry of Space Curves}

\begin{document}

\maketitle

\noindent
{\small \bf Abstract. \rm
The {\em focal curve} of an immersed smooth curve 
$\g:\theta\mapsto \g(\theta)$, in Euclidean space $\R^{m+1}$, consists
of the centres of its osculating hyperspheres. This curve 
may be parametrised in terms of the Frenet frame of $\g$ 
(${\bf t},{\bf n}_1, \ldots,{\bf n}_m$), 
as $C_\g(\theta)=(\g+c_1{\bf n}_1+c_2{\bf n}_2+\cdots+c_m{\bf n}_m)(\theta)$, where 
the coefficients $c_1,\ldots,c_{m-1}$ are smooth functions that we 
call the {\em focal curvatures} of $\g$. We discovered a remarkable formula 
relating the Euclidean curvatures $\k_i$, $i=1,\ldots,m$, of $\g$ with its focal
curvatures. We show that the focal curvatures satisfy a system of Frenet 
equations (not vectorial, but scalar!). We use the properties of the focal curvatures 
in order to give, for $l=1,\ldots,m$, necessary and sufficient conditions for the 
radius of the osculating $l$-dimensional sphere to be critical. We also give necessary
and sufficient conditions  for a point of $\g$ to be a vertex. Finally, 
we show explicitly the relations of the Frenet frame and the Euclidean curvatures of 
$\g$ with the Frenet frame and the Euclidean curvatures of its focal curve $C_\g$.


}
\bigskip


{\bf Introduction}
\medskip

The differential geometry of space 
curves is a classical subject which usually relates geometrical 
intuition with analysis and topology. 
Last years, the ideas and techniques 
of singularity theory of wave fronts and caustics 
(\cite{avg}, \cite{arnoldcw}), revealed to be 
a powerful tool to discover new theorems on the differential geometry 
of curves and surfaces (c.f. \cite{Arnoldlc}-\cite{Arnoldlst},
\cite{Kazarian}, \cite{OT}, 
\cite{Shch}-\cite{Uribe4vtst}). 

The {\em focal surface} or {\em caustic} of a curve $\g$ in Euclidean 
$3$-space is the envelope of the normal planes of $\g$. 
The study of the focal surface of a 
curve can provide useful geometric information about that curve and 
vice versa. Darboux had found how to determine the {\em evolutes} 
of a curve $\g$, that is, the curves whose tangents are normals of 
$\g$. Moreover, he had shown that {\em the focal surface of $\g$ is 
foliated by the evolutes, and all of them lie on 
the focal surface}, see \cite{Darboux}.

The focal surface of $\g$ is singular along a curve $C_\g$ (it has a 
cuspidal edge along $C_\g$) which is called the {\em focal curve} of $\g$ 
(in \cite{Blaschke-L}, it is called the {\em evolute of second type} of $\g$). 
{\em The osculating planes of $C_\g$ are the normal planes of $\g$, and 
the points of $C_\g$ are the centres of the osculating spheres of $\g$}, 
see \cite{Blaschke-L}.

In this paper, we study the geometry of the focal surface, 
focusing on the properties of the focal curve $C_\g$. Using these 
properties, we formulate and prove new results for curves in Euclidean 
$n$-space for arbitrary $n\geq 2$. 

Let $\g: \R \rightarrow \R^{m+1}$ be a smooth curve (a source of light). 
The {\em caustic} of $\g$ (defined as the envelope of the normal lines of $\g$) 
is a singular and stratified hypersurface.  
The {\em focal curve} of $\g$, $C_\g$, is defined as the singular stratum of 
dimension $1$ of the caustic and {\em it consists
of the centres of the osculating hyperspheres of $\g$}. 
Since the centre of any hypersphere tangent to $\g$ at a point 
lies on the normal plane to $\g$ at that point, the focal curve 
of $\g$ may be parametrised using the Frenet frame 
(${\bf t},{\bf n}_1, \ldots,{\bf n}_m$) of $\g$ as follows: 
$C_\g(\theta)=(\g+c_1{\bf n}_1+c_2{\bf n}_2+\cdots+c_m{\bf n}_m)(\theta)$, 
where the coefficients $c_1,\ldots,c_{m-1}$ are smooth functions 
that we call the {\em focal curvatures} of $\g$.  

The Euclidean curvatures of $\g$, $\k_1,\k_2,\ldots,\k_m$, form a system 
of $m$ functions which determine the curve $\g$ up to translation and rotation. 
Let us denote with a prime the derivation with respect to the arc-length parameter.  
We prove that the following formula holds
(Theorem 2): 
$$\k_i={c_1c_1'+c_2c_2'+\cdots+c_{i-1}c_{i-1}' \over c_{i-1}c_i},\ \ 
\mathit{for}\ i\geq 2,$$
showing that the focal curvatures also determine the curve up to 
translation and rotation.

In Theorem 1, we show that {\em the focal curvatures of $\g$ satisfy a system of 
Frenet equations} (not vectorial, but scalar equations and with the same Frenet 
matrix of $\g$ !).  
\medskip

For $k=1,\ldots, m-1$, we give necessary and sufficient conditions, in terms of 
the focal curvatures, for which the radius of the $k$-dimensional osculating 
sphere of a generic curve in $\R^{m+1}$ be critical (Theorem 4).
\medskip

We prove that: {\em A point of $\g$ is a vertex {\rm (that is, a point at which the order 
of contact of $\g$ with its osculating hypersphere is higher than the usual one)} 
if and only if $c_m'+c_{m-1}\k_m=0$ at\ that\ point} (Theorem 3). So, in terms of the focal 
curvatures, the equation characterising the curves lying on a hypersphere 
in $\R^{m+1}$ is very simple: $c_m'+c_{m-1}\k_m\equiv 0$.
\medskip

In Theorem 5, we show explicitly that the Frenet frame 
of the focal curve $C_\g$ consist (up to signs) of the same vectors that 
the Frenet frame of $\g$ but the order of the vectors is inversed. Moreover, 
the Euclidean curvatures $K_1,\ldots,K_m$ of the focal curve $C_\g$ are 
related to those of $\g$ by 
$${K_1 \over |\k_m|}={K_2 \over \k_{m-1}}
=\cdots ={|K_m| \over \k_1}= {1 \over |c_m'+c_{m-1}\k_m|}.$$ 
These relations, together with the stratification 
of the caustic described in \S2, provide a partial solution to the {\em inverse problem}: 
given the caustic, reconstruct the source of light. 
 \medskip

In \S 0, we define the order of contact of a curve with a submanifold of 
$\R^n$ and we recall some basic notions and results on the differential geometry 
of space curves. In \S 1, we state the results of the paper. 
In \S2, we use the techniques of singularity theory (in symplectic geometry) 
to study the geometry and the natural stratification of the {\em focal set} 
of a curve $\g$ in Euclidean $n$-space (the codimension 1 strata being the focal curve
of $\g$). In \S3, we prove our results. 
\medskip


\noindent
{\bf \S 0. Preliminary Definitions and Remarks}
\medskip

In order to give the definition of osculating $k$-spheres of a curve (at a 
point of it) 
we need to introduce the following definition: 
\medskip 

\noindent
{\bf Definition}.~~Let $M$ be a $d$-dimensional submanifold of
$\R^n$, considered as a complete 
intersection: $M=\{ x \in \R^n : g_1(x)=\cdots =g_{n-d}(x)=0 \}$.
We say that a (regularly parametrised) smooth curve
$\g : \theta \mapsto \g(\theta)\in \R^n$ has {\em $k$-point contact}
with $M$ or that their {\em order of contact} is
$k$, at a point $\g(\theta_0)$, if at $\theta=\theta_0$ each function 
$g_1\circ \g,\ldots ,g_{n-d}\circ \g$ has a zero of multiplicity 
at least $k$ and at least one of them has a zero of multiplicity $k$.
\medskip 

\noindent
{\sl Remark}.~~To make this definition more invariant, 
one could denote the image of $\g$ by $\G$ and then write that the 
{\em order of contact} at a point is the minimum of 
the multiplicities of zero among the functions of
the form $g_{|\Gamma}:\Gamma\to\R$ at that point, 
where $g$ belongs to the generating
ideal of $M$ and we assume that $0$ is a regular value of $g$.
\medskip 

In this paper, $M$ will be an
affine subspace or a sphere of dimension $d$. 
\medskip 

\noindent
{\sl Remark}.~~Do not confuse our order of contact with the order 
of tangency: two perpendicular lines in the plane have order of contact $1$ 
at the point of intersection, but the order of tangency is $0$. 
\medskip 


\noindent
{\bf Example 1}.~~A smooth curve in Euclidean (or affine) space $\R^n$ has 
$2$-point contact with its tangent line (at the point of tangency) 
for the generic points of the curve. The plane curve $y=x^3$ has $3$-point contact
  with the line $y=0$, at the origin: the equation $x^3=0$ has a root of
  multiplicity $3$.
\medskip 

\noindent
{\bf Conventions}: Write $n=m+1$. 
In the sequel $\R^{m+1}$ denotes a Euclidean space,
$\theta$ denotes any regular parameter of the curve and $s$
denotes the arc length parameter.  
A parametrised curve $\g=\g(\theta)$ in $\R^{m+1}$ is said to be {\em good} if its 
derivatives of order $1,\dots,m,$ are linearly independent at any point. 
{\em A generic curve is good}. We will consider only good curves.  
\medskip 

The {\em osculating $k$-plane} of a curve at a point is the affine subspace 
spanned by the first $k$ derivatives of the curve 
at that point. {\em A curve has at least ($k+1$)-point contact 
with its osculating $k$-plane at the point of osculation}. 
For $k=m$ we will simply write {\em osculating hyperplane}.

Given a point of a generic smoothly immersed curve in $\R^{m+1}$, 
the sequence consisting of that point and of the osculating 
$k$-planes, $k=1,\ldots,m$, form a complete flag, which is called 
the {\em osculating flag} of the curve at that point. 

By convention, the $k$-dimensional affine subspaces of the Euclidean
space $\R^{m+1}$ will be also considered as $k$-dimensional spheres of
infinite radius. 
\medskip 

\noindent
{\bf Definition}.~~For $k=1,\dots,m,$ a {\em $k$-osculating sphere} at a 
point of a curve in the Euclidean space $\R^{m+1}$ is a $k$-dimensional 
sphere having at least $(k+2)$-point contact with the curve at that point. 
For $k=m$ we will simply write {\em osculating hypersphere}.
\medskip 

\noindent
{\bf Example 2}.~~A generic plane curve and its osculating circle have $3$-point 
contact at an ordinary point of the curve. 
\medskip 

\noindent
{\sl Remark}.~~For $1\leq l<m$, the osculating $l$-sphere at
a point of a curve in $\R^{m+1}$ is the intersection of the osculating 
hypersphere with the osculating ($l+1$)-plane at that point. 
\medskip

\noindent
{\bf Curvature, Frenet frame and higher order curvatures}. \ 
For a curve $\g$ in $\R^3$ parametrised by 
arc-length (from a fixed point) the tangent vector ${\bf t}(s)=\g'(s)$ is unitary 
and it is orthogonal to ${\bf t}'(s)=\g(s)''$. If $\g(s)''\neq 0$ these vectors span the (unique) 
osculating plane of $\g$ at $s$. Write ${\bf t}'(s)=\k_1(s){\bf n}_1(s)$, where ${\bf n}_1(s)$ is 
the unit vector orthogonal to ${\bf t}(s)$ such that the coefficient $\k_1(s)$, called 
the {\em curvature of $\g$ at $s$}, is positive. 
The radius of the osculating circle of $\g$ at $s$ is given by $R_1(s)=1/\k_1(s)$ 
and it is called the {\em radius of curvature of $\g$ at $s$}. 

Assume that $\R^3$ is oriented and take the unit vector ${\bf n}_2(s)$ such that 
the basis ${\bf t}(s)$, ${\bf n}_1(s)$, ${\bf n}_2(s)$, called {\em Frenet frame}, 
is positive (right-handed), that is ${\bf n}_2={\bf t}\times{\bf n}_1$. One easily 
proves that there is a number $\k_2=\k_2(s)$, called the {\em torsion} or {\em second curvature} 
of $\g$ at $s$, such that  ${\bf n}'_2=-\k_2{\bf n}_1$. It is the speed of rotation of 
the vector ${\bf n}_2$. For any good curve we have the following formulas: 
$${\bf t}'=\k_1{\bf n}_1, \ \ \ \ {\bf n}'_1=-\k_1{\bf t}+\k_2{\bf n}_2, \ \ \ \ 
{\bf n}'_2=-\k_2{\bf n}_1,$$
which are called {\em Frenet equations} of the curve $\g$. 
\medskip

Consider a good curve $\g$ in the oriented space $\R^{m+1}$, that is, the vectors 
$\g'(s),\ldots,\g^{(m)}(s)$ are linearly independent for any $s$. Apply Gram-Schmidt 
process to these vectors to obtain the orthonormal system 
${\bf t}(s), {\bf n}_1(s),\ldots, {\bf n}_{m-1}(s)$. Let ${\bf n}_m(s)$ be the 
(unique) vector such that the basis ${\bf t}(s), {\bf n}_1(s),\ldots, {\bf n}_m(s)$, 
called {\em Frenet frame} of $\g$ at $s$, is 
orthonormal and positive. The derivatives of the Frenet frame vectors are given by 
the so called {\em system of Frenet equations} of $\g$: 
$$
\left(
\begin{array}{c}
{\bf t}' \\
{\bf n}_1' \\
{\bf n}_2' \\ 
{\bf n}_3' \\
\vdots \\
{\bf n}_{m-2}' \\
{\bf n}_{m-1}' \\
{\bf n}_m'
\end{array}
\right)=
\left(
\begin{array}{ccccccc}
\ 0 & \ \k_1  & \ 0 & \cdots & 0 & 0 & 0 \\
-\k_1 & \ 0  & \  \k_2 & \cdots & 0 & 0 & 0  \\ 
\ 0   &  -\k_2 & \ 0   & \ & \ & \ & \vdots \\
\ 0   & \ 0  & -\k_3   & \ & \  & \ & \vdots \\
\ \vdots  & \ 0  & \ &  \  & \  & \ & \ \\
\ &  \ & \   & \   & 0 & \k_{m-1} & 0 \\
\ \vdots   &   \ & \   & \   & -\k_{m-1} & 0 & \k_m \\
\ 0     &  \  0  & \  & \cdots & 0 & -\k_m & 0
\end{array}
\right)
\left(
\begin{array}{c}
{\bf t} \\
{\bf n}_1 \\
{\bf n}_2 \\ 
{\bf n}_3 \\
\vdots \\
{\bf n}_{m-2} \\
{\bf n}_{m-1} \\
{\bf n}_m
\end{array}
\right)_. 
$$
The functions $\k_1=\k_1(s),\ldots, \k_m=\k_m(s)$ are called {\em Euclidean curvatures} 
of the curve and are defined only for the good curves. Note that the $l$-th Euclidean curvature
$\k_l$ gives the speed of rotation of the osculating $l$-plane around the osculating 
($l-1$)-plane, with respect to the variation of the arc-length parameter (one can found 
other geometric interpretations of the Euclidean curvatures). 
The curvatures 
$\k_1,\ldots, \k_{m-1}$ of any good curve are strictly positive, while  
$\k_m$ can take any real value. 
\medskip

A point of a smooth curve in $\R^{m+1}$ for which the derivative 
of the curve of order $m+1$ belongs to the osculating hyperplane 
is said to be a {\em flattening}. {\em At a flattening the last Euclidean 
curvature $\k_m$ vanishes and the curve has at least $(m+2)$-point 
contact with its osculating hyperplane at that point}. 
\medskip

\noindent
{\sl Remark about flattenings}.~~At a flattening of a 
generic curve the osculating hypersphere is unique and it coincides 
with the osculating hyperplane. In this case, the centre of the
osculating hypersphere is not defined and we will say that ``it is at
infinity''. If at a point the order of contact of $\g$ with its
osculating sphere of codimension $2$, $S^{m-1}$, is greater than
the usual one, then the point is a non generic flattening. In this case, 
all hyperspheres containing $S^{m-1}$ are osculating, i.e. the centre of 
the osculating hypersphere is not uniquely defined.
\medskip

\noindent
{\bf Example}.~~
These conditions (non satisfied for any point of 
a generic curve) are however satisfied by the flattenings of a
generic spherical curve (that is, a generic curve among the curves lying on a 
hypersphere). 
\medskip

For these reasons we will assume that our curves are good and have no flattening, 
unless we consider (explicitly) spherical curves. 
\medskip


\noindent
{\bf \S 1. Statement of Results}
\medskip

\noindent
{\bf Definition}.~~The curve $C_\g: \theta \mapsto C_\g(\theta)\in \R^{m+1}$
consisting of the centres of the osculating hyperspheres of a good curve 
(without its flattenings) $\g: \theta \mapsto \g(\theta)\in \R^{m+1}$ is called the 
{\em parametrised focal curve} of $\g$. 
\medskip

\noindent
{\sl Remark}.~~In geometrical optics, a curve $\g$ in 
Euclidean $3$-space can be considered as 
a source of light. The envelope of all light rays normal to $\g$ is the
{\em focal surface} or {\em caustic} of $\g$. The light intensity is 
much more concentrated on the caustic than in all other points 
of the space. Moreover, the caustic itself is more illuminated along its 
cuspidal edge, which is the focal curve of $\g$. 
\medskip

Consider a good curve $\g: \R \rightarrow \R^{m+1}$. 
Write $\k_1,\k_2,\ldots,\k_m$ for its Euclidean curvatures and 
${\bf t},{\bf n}_1, \ldots,{\bf n}_m$ for its Frenet frame. 
The hyperplane normal to $\g$ at a point consists of the set of centres of all 
hyperspheres tangent to $\g$ at that point. Hence the centre of the 
osculating hypersphere at that point lies in such normal hyperplane. Therefore 
(denoting $C_\g(\theta)$ by $C_\g$, $\g(\theta)$ by $\g$ and so on,\ldots) we can write 

$$C_\g=\g+c_1{\bf n}_1+c_2{\bf n}_2+\cdots+c_m{\bf n}_m,$$
where the coefficients $c_1,\ldots,c_{m-1}$ are smooth functions of the parameter of the
curve $\g$. 
\medskip

\noindent
{\bf Definition}.~~The coefficient $c_i$, $i=1,\ldots,m$, is called the $i^{th}$
{\em focal curvature} of $\g$. 
\medskip

\noindent
{\sl Remark}.~~{\em The first focal curvature $c_1$ never vanishes: $c_1=1/\k_1$}. 
\medskip

The Frenet equations of a curve in ($m+1$)-Euclidean space is a system of 
$m+1$ vectorial equations involving the unit vectors of the Frenet 
frame and their derivatives. The following theorem shows that the focal curvatures 
of that curve satisfy a system of {\em scalar Frenet equations} which ``is obtained 
from the usual Frenet equations by replacing the $i^{th}$ normal vector of the 
Frenet frame by the $i^{th}$ focal curvature''. 
\medskip

\noindent
{\bf Theorem 1}.~~{\em The focal curvatures of 
a curve lying on a hypersphere
$\g: \R \rightarrow \sph^n\subset \R^{m+1}$, parametrised by arc length $s$, satisfy the
  following ``scalar Frenet equations'': 

$$
\left(
\begin{array}{c}
1 \\
c_1' \\
c_2' \\
c_3' \\
\vdots \\
c_{m-2}' \\
c_{m-1}' \\
c_m'
\end{array}
\right)=
\left(
\begin{array}{ccccccc}
\ 0 & \ \k_1  & \ 0 & \cdots & 0 & 0 & 0 \\
-\k_1 & \ 0  & \  \k_2 & \cdots & 0 & 0 & 0  \\ 
\ 0   &  -\k_2 & \ 0   & \ & \ & \ & \vdots \\
\ 0   & \ 0  & -\k_3   & \ & \  & \ & \vdots \\
\ \vdots  & \ 0  & \ &  \  & \  & \ & \ \\
\ &  \ & \   & \   & 0 & \k_{m-1} & 0 \\
\ \vdots   &   \ & \   & \   & -\k_{m-1} & 0 & \k_m \\
\ 0     &  \  0  & \  & \cdots & 0 & -\k_m & 0
\end{array}
\right)
\left(
\begin{array}{c}
0 \\
c_1 \\
c_2 \\ 
c_3 \\
\vdots \\
c_{m-2} \\
c_{m-1} \\
c_m
\end{array}
\right)_. 
$$}

\noindent
{\sl Remark}.~~If the curve is not spherical then
the correcting term $-\frac{(R_m^2)'}{2c_m}$ must be added to the last component of the
left hand side vector to obtain $c_m'-\frac{(R_m^2)'}{2c_m}$, for $c_m\neq 0$. 
\medskip

\noindent
{\bf Theorem 2}.~~{\em The Euclidean curvatures of a good curve $\g$ 
(with $\k_m\neq 0$) in $\R^{m+1}$, parametrised by arc length, 
are given in terms of the focal curvatures of $\g$ by the
formula: 
$$\k_i={c_1c_1'+c_2c_2'+\cdots+c_{i-1}c_{i-1}' \over c_{i-1}c_i},\ \ 
\mathit{for}\ i\geq 2.$$}

\noindent
{\sl Remark}.~~For a generic curve, the focal curvatures 
$c_i$ or $c_{i-1}$ can vanish at isolated points. At these points the function 
$c_1c_1'+c_2c_2'+\cdots+c_{i-1}c_{i-1}'$ also vanishes, and the
corresponding value of the Euclidean curvature $\k_i$ may be obtained by
l'H{\^o}pital rule.  
\medskip

\noindent
{\bf Definition}.~~A {\em vertex} of a curve in $\R^n$ 
is a point at which the curve has at least $(n+2)$-point contact with its 
osculating hypersphere.
\medskip 

\noindent
{\bf Example 3}.~~{\em The vertices of a curve in Euclidean plane $\R^2$ are the points 
at which the curvature is critical}. For instance, a non-circular ellipse 
has $4$ vertices: They are the points at which the ellipse intersects its
principal axes. 
\medskip 

The interest on the vertices of curves came, for instance, from geometrical 
optics (c.f. Huygens) and from the {\em geometry in the large}. Namely 
the classical $4$-vertex theorem states that {\em a smooth closed convex 
plane curve has at least $4$ different vertices}, \cite{Mukho}. Besides several important works 
generalising this theorem (c.f. \cite{Kneser,Mohrmann,Barner,Segre,Carmen1,Sedykh4f}), 
the recent progress in symplectic geometry and 
singularity theory have revived the interest on the study of vertices together with 
the different variants of its definition
(c.f. \cite{Shch,Ghys,Kazarian,OT}, \cite{Uribe2}-\cite{Uribe4vtst}). 
Here we are mainly concerned with local properties of vertices.

The next theorem (implicitly contained in \cite{Carmen4}) provides necessary and sufficient 
conditions for a point to be a vertex. 
\medskip

\noindent
{\bf Theorem 3}.~~{\em A non-flattening point of a good curve parametrised by arc
  length in $\R^{m+1}$, $m>1$, is a vertex if and
  only if $$c_m'+c_{m-1}\k_m=0 \ \ \mathit{at\ that\ point.}$$}

\noindent
{\bf Corollary 1}.~~{\em A good curve parametrised by arc length 
in Euclidean space $\R^{m+1}$, $m>1$, lies on a hypersphere if and
  only if $$c_m'+c_{m-1}\k_m\equiv 0.$$}

\noindent
{\bf Example 4}.~~For curves in Euclidean $3$-space, Corollary 1 provides the following 
classical result on spherical curves (see for instance \cite{Blaschke}): 
\medskip

{\em A smoothly immersed curve of $\R^3$, with curvature $\k$ and torsion 
$\tau$ both nowhere zero, lies on a sphere if and only if 
$$c_2'+c_1\tau \equiv 0, \mbox{ i.e. if and only if } 
\left({R_1' \over \tau}\right)'+R_1\tau \equiv 0,$$
where derivation is taken with respect to the arc length of the curve and 
$R_1=1/\k$, is the radius of curvature.} 
\medskip

Unfortunately, I have found a small mistake in the 
beautiful Hilbert--Cohn Vossen's book, \cite{Hilbert} :
\medskip

{\em A curve of $\R^3$ lies on a sphere if and only if 
$$R_1^2+({R_1'})^2{1 \over \tau^2} = \mathrm{const}.
\eqno(W)$$}

Of course, a curve lying on a sphere satisfies condition (W), 
which means that the radius of the osculating sphere is constant. 
However, the number of non-spherical curves satisfying condition 
(W) is infinite: {\em If a curve with nowhere vanishing torsion has 
constant curvature $\k\neq 0$ then the radius of its osculating sphere 
is constant and equal to $R=1/\k$.} This follows from condition (W). 
One example is the circular helix $t\mapsto (\cos t, \sin t, t)$). 
The above statement becomes true if one suppose the genericity 
condition $R_1'\neq 0$. 
\medskip

{\em The radius of the osculating hypersphere of a curve in $\R^{m+1}$ 
is critical at each vertex of that curve; the converse statement 
is not always true for $m>1$} (see \cite{Uribe1}, \cite{Uribe4vtst}):
There are examples of curves having points for which the radius of the 
osculating hypersphere is critical, but which are not vertices. The geometric 
meaning of such points becomes clear from Theorem 5, below.

The following two theorems give necessary 
and sufficient conditions for the radius of the osculating sphere 
of dimension $l \leq m$ to be critical. 
\medskip

\noindent
{\bf Theorem 4}.~~{\em For $1\leq l<m$, the radius of the
  osculating $l$-sphere of a generic curve in $\R^{m+1}$
is critical if and only if either 

\centerline{$c_l=0$\ ~ or\ ~ $c_{l+1}=0$.}

\noindent
Moreover, $c_1$ never vanishes.} 
\medskip

\noindent
{\sl Remark}.~~{\em At a point of a curve $\g$, the first $l$ 
focal curvatures $c_1,\ldots,c_l$ are the coordinates (with respect to the Frenet frame) 
of the centre of the $l$-dimensional osculating sphere of $\g$ at 
that point.} Therefore the curve $\g_l$ described by the centre of the $l$-dimensional 
osculating sphere is parametrised by: 
$$\g_l=\g+c_1{\bf n}_1+c_2{\bf n}_2+\cdots+c_l{\bf n}_l.$$
Of course, $\g_m=C_\g$. Theorem 4 implies for instance that {\em the curves 
$\g_1$ and $\g_2$ intersect at least twice and the curve $\g_l$ intersects 
either $\g_{l-1}$ or $\g_{l+1}$, at least at two points, $1<l<m$.}
\medskip

\noindent
{\bf Corollary 2}.~~{\em If the $l^{th}$ focal curvature $c_l$ vanishes
  at a point, then both the radii of the osculating spheres 
of dimensions $l-1$ and $l$ are critical at that point.} 
\medskip

\noindent
{\bf Remark}.~~The $m^{th}$ focal curvature $c_m$ at a point 
of a smooth curve in $\R^{m+1}$ is the signed distance between the osculating 
hyperplane and the centre of the osculating hypersphere at that point. 
\medskip

\noindent
{\bf Definition}.~~A point of a curve is said to be 
a {\em pseudo-vertex} of that curve if the centre of the osculating 
hypersphere at that point lies in the osculating hyperplane at 
that point (that is, if $c_m=0$). 
\medskip

\noindent
{\bf Corollary 3}.~~{\em A generic closed curve in $\R^3$ has at least 
two vertices or two pseudo-vertices}. 
\medskip

\noindent
{\bf Corollary 4}.~~{\em At a pseudo-vertex of a smooth curve in 
$\R^{m+1}$, $m>1$, both the radius of the osculating hypersphere and the 
radius of the osculating $(m-1)$-sphere are critical}. 
\medskip

\noindent
{\bf Proposition 0}.~~{\em The radius of the osculating hypersphere
 at a point  of a good curve in $\R^{m+1}$, $m>1$, is critical 
if and only if such point is either a vertex or a pseudo-vertex.} 
\medskip

A point of a generic smooth curve at which the last Euclidean curvature vanish, $\k_m=0$, 
is a {\em flattening} of the curve (see our Remark about flattenings above). The 
following statement is a consequence of Proposition 0. 
\medskip

\noindent
{\bf Corollary 5}.~~{\em Write $V$, $F$ and $P$ for the  number of
  vertices, flattenings and pseudo-vertices of a generic closed curve
  smoothly immersed in $\R^{m+1}$. The following inequalities hold: 
$$V+P\geq F \mbox{~  and ~} V+P \geq 2.$$}

We reformulate Proposition 0 (and we will prove it, in \S 3) in terms of the focal 
curvatures $c_m$ and $c_{m-1}$: 
\medskip

\noindent
{\bf Proposition \~{0}}.~~{\em The radius of the osculating hypersphere
  of a good curve in $\R^{m+1}$, $m>1$, parametrised by arc length, 
is critical at a point 
if and only if either $c_m=0$ or $c_m'+c_{m-1}\k_m=0$ at that point.}
\medskip

After I have sent this paper to V.D. Sedykh, he communicated to me 
that he had discovered independently Proposition 0 and Corollary 5, but he had not 
published them and he urged me to publish all results of this paper. 
\medskip

\noindent
{\sl Remark}.~~By definition, the first $m-1$ Euclidean curvatures of a generic
curve $\g : \R \rightarrow \R^{m+1}$ are positive everywhere, while 
the last one, $\k_m$, can take any real value. 
The sign of the last Euclidean curvature at a non-flattening point of a curve 
is defined only when 
the orientation on the ambient space $\R^{m+1}$ is fixed : 
$\k_m$ is positive (negative) at the points of the curve where 
the derivatives of order $1,\ldots, m+1$ form a positive (negative, resp.) 
basis of $\R^{m+1}$. 
\medskip

\noindent
{\sl Remark}.~~Consider a curve $\g : \R \rightarrow \R^{m+1}$ in the 
oriented Euclidean space $\R^{m+1}$. 
{\em If the number $m>0$ is of the form $4k$ or $4k+1$, with 
$k\in \N$, then sign of the last Euclidean curvature of $\g$ at a non-flattening point 
depends on the orientation of the curve}. That is, the last Euclidean curvature 
of a curve at a non-flattening point is a function whose sign depends not only 
on the point of the curve but also on the orientation of the curve 
given by the parametrisation. 
\medskip

\noindent
{\sl Proof}.~~Let $\g : \R \rightarrow \R^{m+1}$ be a generic curve in $\R^{m+1}$, 
such that $\g(0)$ is not a flattening. 
Write $\tau(t) = -t$ and 
consider the parametrisation in the opposite direction 
$\tilde{\g} =\g \circ \tau: t \mapsto \g(-t)$. The derivative of order 
$r$ of $\tilde{\g}$ at $t=0$ is $\tilde{\g}^{(r)}(0)=\g^{(r)}(0)\cdot (-1)^r$. 
So the derivatives of odd order of $\g$ and $\tilde{\g}$ at $t=0$ have opposite 
directions while the derivatives of even order of $\g$ and $\tilde{\g}$ at 
$t=0$ coincide. Therefore the basis obtained from 
the derivatives of order $1,\ldots, m+1$ of $\tilde{\g}$ at $t=0$ 
and the basis obtained from the derivatives of order $1,\ldots, m+1$ of $\g$ at 
$t=0$ give different orientations of $\R^{m+1}$ if and only if 
the cardinality of the set $\{ r\in \N : \mbox{ $r$ is odd and $r\leq m+1$ }\}$ 
is odd, i.e. 
if and only if the number $m>0$ is of the form $4k$ or $4k+1$, with 
$k\in \N$. 
\medskip

\noindent
{\bf Theorem 5}.~~{\em Let $\g : s \mapsto \g(s)\in \R^{m+1}$ be a good curve 
without its flattenings. 
Write $\k_1,\ldots,\k_m$ for its Euclidean curvatures and 
$\{ {\bf t},{\bf n}_1, \ldots,{\bf n}_m\}$ for its Frenet frame. 
For each non-vertex $\g(s)$ of $\g$, write $\e(s)$ for the sign of 
$(c_m'+c_{m-1}\k_m)(s)$ and $\d_k(s)$ for the sign of $(-1)^k\e(s)\k_m(s)$, 
$k=1,\ldots,m$. For any non-vertex of $\g$ the following holds: 
\medskip

a) The Frenet frame 
$\{ {\bf T},{\bf N}_1, \ldots,{\bf N}_m\}$ 
of $C_\g$ at $C_\g(s)$ is well-defined and its vectors are given by 
${\bf T}=\e {\bf n}_m$, ${\bf N}_k=\d_k {\bf n}_{m-k}$, for 
$k=1,\ldots,m-1$, and ${\bf N}_m=\pm {\bf t}$, 
the sign in $\pm {\bf t}$ is chosen in order to obtain a positive 
basis. 
\medskip

b) The Euclidean curvatures $K_1,\ldots,K_m$ of 
the parametrised focal curve of $\g$, $C_\g :s \mapsto C_\g(s)$,
are related to those of $\g$ by :
$${K_1 \over |\k_m|}={K_2 \over \k_{m-1}}
=\cdots ={|K_m| \over \k_1}= {1 \over |c_m'+c_{m-1}\k_m|},$$
the sign of $K_m$ is equal to $\d_m$ times 
the sign chosen in $\pm {\bf t}$. 
\medskip

That is, the Frenet matrix of $C_\g$ at $C_\g(s)$ is 

$$
{1 \over |c_m'+c_{m-1}\k_m|}
\left(
\begin{array}{ccccccc}
\ 0 & |\k_m|  & \ 0 & \cdots & 0 & 0 & 0 \\
-|\k_m| & \ 0  & \  \k_{m-1} & \cdots & 0 & 0 & 0  \\ 
\ 0   &  -\k_{m-1} & \ 0   & \ & \ & \ & \vdots \\
\ 0   &  0  & -\k_{m-2}   & \ & \  & \ & \vdots \\
\ \vdots  & \ 0  & \ &  \  & \  & \ & \ \\
\ &  \ & \   & \   & 0 & \k_2 & 0 \\
\ \vdots   &   \ & \   & \   & -\k_2 & 0 & \mp \d_m \k_1 \\
\ 0     &  \  0  & \  & \cdots & 0 & \pm \d_m \k_1 & 0
\end{array}
\right)_. 
$$}
\medskip

\noindent
{\bf Application to self-congruent curves}. \ 
A curve of $\R^{m+1}$ is said to be {\em self-congruent} if for any two 
points $a$ and $b$ of it, there is a preserving orientation orthogonal 
transformation of $\R^{m+1}$ sending the curve to itself and sending 
$a$ to $b$. One can prove that {\em the class of self-congruent curves coincides 
with the class of curves whose Euclidean curvatures are constant}. 

The focal curvatures of these curves are therefore constant and the scalar Frenet 
equations imply that 
$$c_{2l}=0 \ \ \mbox{and } \ \ c_{2l+1}=\prod_{j=0}^l\left(\frac{k_{2j}}{k_{2j+1}}\right),$$ 
where the convention $\k_0=1$ is used, and the subindices $2l$ and $2l+1$ are taken 
over all values of $l$ for which $2\leq 2l\leq m$ and $1\leq 2l+1 \leq m$, respectively. 
\medskip

\noindent 
{\bf Proposition}.~~
{\em For any $l\in \N$ such that $0<2l\leq m$, 
the following holds: At any point of a self-congruent curve of $\R^{m+1}$ the 
centre of the osculating $2l$-sphere lies in the osculating $2l$-plane}. 
\medskip

\noindent
{\em Proof}. \ This follows from the above equalities $c_{2l}=0$. 
\bigskip


\noindent
{\bf \S 2. Study of the Focal Set (caustic) of a Curve}
\medskip

The {\em focal set} or {\em caustic} of a submanifold of positive
codimension in Euclidean space $\R^{m+1}$ (for instance, of 
a curve in $\R^3$) is defined as the envelope of the family of 
normal lines to the submanifold. 
\medskip

\noindent
{\sl Remark}.~~Similarly to geometrical optics in Euclidean $3$-space, 
a submanifold of positive codimension in Euclidean space $\R^{m+1}$ 
may be considered as a source of light (or as an initial wave front). 
The normal lines to this source submanifold are called 
{\em normal light rays} and its focal set (on which the light 
intensity is much more concentrated than in the other points 
of the space) is called the {\em caustic} of that submanifold. 

We will study the focal set of a generic curve 
$\g:\R \rightarrow \R^{m+1}$. 

The hyperplane normal to $\g$ at a point is the union of all lines 
normal to $\g$ at that point. The envelope of all hyperplanes normal 
to $\g$ is thus a component of the focal set that we call the 
{\em main component} (the other component 
is the curve $\g$ itself, but we will not consider it). 

The normal hyperplanes of a curve at two neighbouring points intersect 
along an affine subspace of codimension 2 which approaches a limiting 
position as the points move into coincidence. The affine subspace that 
assumes this limiting position is called the {\em $2$-codimensional 
focal subspace} of the curve at the point under consideration. 

When the point moves along the curve the $2$-codimensional  
focal subspace generates a hypersurface which, by construction, is the envelope 
of the hyperplanes normal to $\g$, i.e. it is the main component of 
the focal set. 

So the main component of the focal set of a curve is the union 
(in a one-parameter family) of affine subspaces of codimension 2 
(see Claim 3 in subsection 2.2). 
\medskip

\noindent
{\bf Example 5}.~~At a point of a curve in $\R^3$, the $2$-codimensional 
focal subspace is the line through the centre of the osculating circle 
which is parallel to the bi-normal vector. In classical differential 
geometry of curves in Euclidean $3$-space, it is called 
the {\em polar line} (see \cite{Darboux}). 
\medskip

\noindent
{\bf 2.1 The caustic of a family of functions}.~~We 
will use techniques of singularity theory in order to have a 
more detailed study of the focal set. 
\medskip

\noindent
{\bf Definition}.~~The {\em caustic} of a family of functions 
depending smoothly on parameters consists of the parameter values 
for which the corresponding function has a non-Morse critical point.
\medskip

\noindent
{\bf Example 6}.~~Given a generic curve $\g:\R \rightarrow \R^{m+1}$, 
let $F:\R^{m+1} \times \R \rightarrow \R$ be the $(m+1)$-parameter family 
of real functions given by  
$$F(q,\theta)={1 \over 2}\parallel q-\gamma(\theta) \parallel^2.$$
The caustic of the family $F$ is given by the set 
$$\{ q\in \R^{m+1}:\exists \theta \in \R : F_q'(\theta)=0\ \mathrm{and}\ F_q''(\theta)=0\}.$$

\noindent
{\bf Proposition A}.~~{\em The caustic of the family 
$F(q,\theta)={1 \over 2}\parallel q-\gamma(\theta) \parallel^2$ coincides with the 
focal set of the curve $\g:\R \rightarrow \R^{m+1}$}. 
\medskip

\noindent
{\sl Proof}.~~The caustic of $F$ is defined by the pair of equations 
$F_q'(\theta)=0$, $F_q''(\theta)=0$. For each fixed value of $\theta$, the set of 
points $q\in \R^{m+1}$ satisfying the first equation form the hyperplane 
normal to $\g$ at $\g(\theta)$ :
$$F_q'(\theta)= -\langle q-\g(\theta), \g'(\theta) \rangle = 0.$$
The set of points $q\in \R^{m+1}$ satisfying both equations for a fixed $\theta$ 
are thus the stationary points of the normal hyperplane at $\g(\theta)$ under an 
infinitesimal variation of it. They form an affine subspace of 
codimension 2 in $\R^{m+1}$ : 
$$F_q''(\theta)=-\langle q-\g(\theta), \g''(\theta) \rangle + 
\langle \g'(\theta), \g'(\theta) \rangle = 0.$$
Of course this subspace coincides with the $2$-codimensional 
focal plane of the curve at $\g(\theta)$, considered above. 
\medskip

\noindent
{\bf 2.2 The natural stratification of the focal set}.~~The focal set
of a curve $\g:\R \rightarrow \R^{m+1}$ is stratified 
in a natural way. The following claims describe the geometry of such 
stratification for curves without flattenings.
\medskip

Denote by $A_\g^k(\theta)$, $k=1,\ldots,m+2$, the 
set consisting of the centres of all hyperspheres having at least 
$(k+1)$-point contact with $\g$ at $\g(\theta)$. 
\medskip

\noindent
{\bf Claim 1}.~~{\em The set $A_\g^k(\theta)$, $k=1,\ldots,m+1$ is an affine 
subspace of codimension $k$ in $\R^{m+1}$.}
\medskip

\noindent
{\bf Claim 2}.~~{\em The set $A_\g^1(\theta)$ (consisting  of the centres 
of all hyperspheres having at least $2$-point contact with $\g$ at $\g(\theta)$) 
is the hyperplane normal to $\g$ at the point $\g(\theta)$.}
\medskip

\noindent
{\bf Definition}.~~The affine subspace $A_\g^k(\theta)$ is 
called {\em $k$-codimensional focal plane} of $\g$ at $\g(\theta)$. 
\medskip

\noindent
{\bf Corollary} (of claims 1 and 2).~~{\em  The sequence of 
focal subspaces $A_\g^1(\theta)\supset A_\g^2(\theta)\supset  \cdots \supset A_\g^{m+1}(\theta)$ 
defines a complete flag on the hyperplane normal to $\g$ 
at $\g(\theta)$.}
\medskip

\noindent
{\sl Remark}.~~The complete flag
$A_\g^1(\theta)\supset A_\g^2(\theta)\supset  \cdots \supset A_\g^{m+1}(\theta)$ 
defines a natural stratification on the hyperplane normal to $\g$ at 
$\g(\theta)$. This stratification induces a natural 
stratification on the focal set of $\g$. The stratum of dimension $1$ 
being the focal curve of $\g$. The $0$-dimensional stratum consists of 
isolated points at which the focal curve is singular (it has a cusp, 
see Proposition 1 in \S 3). These singular points of the focal curve of 
$\g$ correspond to the vertices of $\g$ (for these points the set 
$A_\g^{m+2}(\theta)$ is not empty). 
\medskip

\noindent
{\bf Claim 3}.~~{\em The focal set of a smooth curve
consists of the centres of all hyperspheres having at least 
$3$-point contact with that curve at a point of it
(i.e. it is the union of all the $2$-codimensional focal planes of 
the curve).} 
\medskip

\noindent
{\bf Proposition B}.~~{\em The complete flag 
$A_\g^1(\theta)\supset A_\g^2(\theta)\supset  \cdots \supset A_\g^{m+1}(\theta)$ 
is the osculating flag of the focal curve of $\g$ at the point $C_\g(\theta)$. 
In particular, the hyperplane normal to $\g$ at $\g(\theta)$ coincides 
with the osculating hyperplane of the focal curve of $\g$ at the 
point $C_\g(\theta)$.} 
\medskip

\noindent
{\bf Lemma 0}.~~{\em A point $q\in \R^{m+1}$ is the centre of a hypersphere 
having $k$-point contact with $\g$ at the point $\g(\theta_0)$ if and 
only if the function 
$F_q(\theta)={1 \over 2}\parallel q-\gamma(\theta) \parallel^2$ 
has a critical point of multiplicity $k-1$ at $\theta_0$:
$$F'_q(\theta_0)= F''_q(\theta_0)= \ldots =F^{(k-1)}_q(\theta_0)= 0 \mbox{ and }
F^k_q(\theta_0)\neq 0.$$} 

\noindent
{\sl Proof}.~~The sphere of radius $r$ with centre at $q$ is defined 
by the equation 
$$g_r(x)={1 \over 2}(\parallel q-x \parallel^2-r^2) = 0.$$
So a point $q$ is the centre of a hypersphere having 
$k$-point contact with $\g$ at the point $\g(\theta_0)$ if and 
only if the function $g_r\circ \g$ has a zero of multiplicity $k$ 
at $\theta=\theta_0$, for some $r$, i.e. if and only if the function 
$F_q(\theta)={1 \over 2}\parallel q-\gamma(\theta) \parallel^2$ 
has a critical point of multiplicity $k-1$ at $\theta_0$. \fin
\medskip

\noindent
{\sl Proof of claims 2 and 3}.~~To prove Claims 2 and 3, use Lemma 0
and repeat the proof of Proposition A. Another proof of Claim 3 follows 
from Example 6, Lemma 0 and Proposition A. \fin
\medskip

\noindent
{\sl Proof of claim 1}.~~
Consider the following system of $(m+2)$ equations
$$
\begin{array}{rcr}
F_q'(\theta)  &  =   & 0 \\
F_q''(\theta)  &  =   & 0 \\ 
         & \vdots &  \\ 
F_q^{(m+2)}(\theta)  &  =   & 0. 
\end{array}
$$

For each fixed value of $\theta$, it can be easily seen that 
the first $k$ equations ---written explicitly--- define an affine subspace 
of codimension $k$ in $\R^{m+1}$ (the cases $k=1,2$, are in the proof 
of Proposition A). So the set $A_\g^k(\theta)$ of centres of all hyperspheres 
having at least $(k+1)$-point contact with $\g$ at $\g(\theta)$ is an 
affine subspace of $\R^{m+1}$. \fin
\medskip

\noindent
{\sl Remark}.~~The (generating) family 
$F(q,\theta)={1 \over 2}\parallel q-\gamma(\theta) \parallel^2$ together with 
Sturm theory can be used to calculate the number of vertices 
of the curve $\g$, see \cite{Uribe4vtst}.
\medskip

\noindent
{\sl Remark ({\small for Singularity Theory Specialists})}.~~In the
setting of the theory of Lagrangian singularities, Lagrangian maps 
and the caustics of Lagrangian maps, the focal set of the curve $\g$ 
is the caustic of the {\em Normal map} associated to $\g$, which is a 
Lagrangian map defined by the generating family $F(q,\theta)$ (for the notions 
of caustic, Lagrangian map, Lagrangian singularity and 
generating family, we refer the reader to \cite{avg} and \cite{arnoldcw}). 
Thus {\em the vertices of a curve in $\R^{m+1}$ correspond to a Lagrangian singularity 
$A_{m+2}$ of the normal map, that is, the focal set has a ``swallowtail'' singularity 
at the centres of the osculating hyperspheres corresponding to the vertices of the curve.} 
\medskip



\noindent
{\bf \S 3. The Proofs of the Results}
\medskip

As we mentioned in the introduction, the ideas and techniques of the theory 
of Lagrangian and Legendrian singularities (singularities of caustics and 
wave fronts) were an important tool for the discovery of the results of this 
paper and also for their initial proofs.
Some of these results would be difficult to discover only using 
Frenet frame theory. 
However, once the results were discovered and proved, 
the author has made an effort in order to present the proofs as short as possible 
and as elementary as possible. The author hopes the proofs will be understandable 
for anyone. 
\medskip

To prove our results we will prove before some lemmas related to the
focal curve. Below, $\theta$ denotes any regular parameter of the curve and $s$
denotes the arc length parameter.
\medskip

\noindent
{\bf Lemma 1}.~~{\em Let 
$\g: \theta \mapsto (\f_1(\theta),\ldots,\f_{m+1}(\theta))$ be a good curve in
$\R^{m+1}$. The velocity vector $q'(\theta)$ of the focal curve of $\g$ at
$\theta$ is proportional to the $m^{th}$-normal vector ${\bf n}_m(\theta)$ of
$\g$.}
\medskip

\noindent
{\sl Proof}.~~Consider the (generating) family of functions 
$F:\R \times \R^{m+1} \rightarrow \R$ defined by 
$$F_q(\theta)={1 \over 2}\parallel q-\gamma(\theta) \parallel^2.$$

Write $g={\g^2 \over 2}$. As in \S 2, use the fact that 
$-F=\g \cdot q -{\g^2 \over 2}- {q^2 \over 2}$ to 
recall that the following system of $m+1$ equations defines the 
{\em focal curve} $q(\theta)$ of $\g$: 

$$
\begin{array}{rcl}
\g'\cdot q(\theta)- g'  &  =   & 0, \\
\g''\cdot q(\theta)- g''  &  =   & 0, \\ 
         & \vdots &  \\ 
\g^{(m+1)}\cdot q(\theta)- g^{(m+1)}  &  =   & 0. 
\end{array}
\eqno(*)$$

Derive each equation with respect to $\theta$ to obtain a
second system of equations: 

$$
\begin{array}{rcl}
\g'\cdot q'(\theta) + \g''\cdot q(\theta)- g'' &  =   & 0, \\
\g''\cdot q'(\theta)+ \g'''\cdot q(\theta)- g'''   &  =   & 0, \\ 
         & \vdots &  \\ 
\g^{(m)}\cdot q'(\theta)+ g^{(m+1)}\cdot q(\theta)- g^{(m+1)}  &  =   & 0, \\
\g^{(m+1)}\cdot q'(\theta)+ g^{(m+2)}\cdot q(\theta)- g^{(m+2)}  &  =   & 0. 
\end{array}
\eqno(**)$$

Combine the $i^{th}$ equation of system $(**)$ with the
$(i+1)^{th}$ equation of system $(*)$, for $i=1,\ldots,m$, to obtain 

$$
\begin{array}{rcl}
\g'\cdot q'(\theta) & =   & 0, \\
\g''\cdot q'(\theta)  &  =   & 0, \\ 
         & \vdots &  \\ 
\g^{(m)}\cdot q'(\theta)  &  =   & 0. 
\end{array}
\eqno(***)$$

This means that the velocity vector $q'(\theta)$ is orthogonal to the
osculating hyperplane of $\g$, i.e. $q'(\theta)$ is proportional to
the $m^{th}$-normal vector ${\bf n}_m$. \fin
\medskip

\noindent
{\bf Proposition 1}.~~{\em A non-flattening point of a good curve in
  $\R^{m+1}$ is a vertex if and only if the velocity vector of the
  focal curve is zero.}
\medskip

\noindent
{\sl Proof}.~~If the point $\g(\theta)$ is a vertex of $\g$, then besides
the system of equations $(*)$ obtained in the proof of Lemma 1, it
also satisfies the equation: 

$$\g^{(m+2)}\cdot q(\theta)- g^{(m+2)}  =  0,$$
which combined with the last equation of system $(**)$ gives the
equation 
$$\g^{(m+1)}\cdot q'(\theta)  = 0.$$

The preceding equation together with the system $(***)$ imply that for
a non-flat vertex $\g(\theta)$ of the curve $\g$ the velocity vector
$q'(\theta)$ of the focal curve is zero.  

Conversely, if a point $\g(\theta_0)$ is not a vertex then the corresponding
point of the focal curve satisfies the relation

$$\g^{(m+2)}(\theta_0)\cdot q(\theta_0)- g^{(m+2)}(\theta_0) \neq  0,$$
which together with the last equation of $(**)$, for $\theta=\theta_0$, imply
that $q'(\theta_0)\neq 0.$ \fin
\medskip

Lemma 1 and Proposition 1 were
also stated in \cite{Carmen4}, where the condition to
the point to be a non-flattening is unfortunately absent. 
Without this condition Proposition 1 does not hold. 
\medskip

\noindent
{\bf Lemma 2}.~~{\em Let $\g: \R \rightarrow \R^{m+1}$ be a good curve with $\k_m\neq 0$.  
The derivative of its parametrised focal curve $C_\g$ with respect 
the arc length $s$ of $\g$ is 
$$C_\g' = (c_m'+c_{m-1}\k_m){\bf n}_m.$$}

\noindent
{\bf Proof of Theorem 1, Proposition 0 and Lemma 2}.~~Consider the 
parametrised focal curve of $\g$:  
$$C_\g(s)=(\g+c_1{\bf n}_1+c_2{\bf n}_2+\cdots+c_m{\bf n}_m)(s).$$ 
Denote $C_\g(\theta)$, $\g(\theta)$ and so on by 
$C_\g$, $\g$, etc. Derive $C_\g$ with respect to the arc length of $\g$ 
and use Frenet equations of $\g$ to obtain: 
$$
\begin{array}{rcl}
C_\g'  & =  & {\bf t}+c_1(-\k_1{\bf t}+\k_2{\bf n}_2)+c_1'{\bf n}_1
+\cdots+c_{m-1}'{\bf n}_{m-1}+ c_m(-\k_m{\bf n}_{m-1})+c_m'{\bf n}_m \\
 \ & = & (1-c_1\k_1){\bf t}+(c_1'-\k_2c_2){\bf n}_1+
(c_2'+c_1\k_2-c_3\k_3){\bf n}_2+\cdots \\
 \ & \ & + (c_i'+c_{i-1}\k_i-c_{i+1}\k_{i+1}){\bf n}_i+\cdots+ 
(c_m'+c_{m-1}\k_m){\bf n}_m. 
\end{array}
$$
By Lemma 1, the first $m-1$ components of $C_\g'$ vanish. 
Consequently
$$C_\g'=(c_m'+c_{m-1}\k_m){\bf n}_m
\eqno(1)$$ 
and the following equalities hold:  
$$
\begin{array}{rcl}
1    & = & \k_1c_1, \\
c_1' & = & \k_2c_2, \\
c_2' & = & -\k_2c_1+\k_3c_3, \\
\vdots \ & \vdots & \ \vdots \\
c_{m-1}' & = & -c_{m-2}\k_{m-1}+c_m\k_m. 
\end{array}
\eqno(2)$$
Equation (1) proves Lemma 2. 
Use the fact that the
radius $R_m$ of the osculating hypersphere satisfies 
$R_m^2=\parallel C_\g -\g \parallel^2$ to obtain 
$$
\begin{array}{rcl}
(R_m^2)' & = & \langle C_\g-\g, C_\g-\g \rangle' \\
\ &   = & 2\langle C_\g'-\g', C_\g-\g \rangle \\
\ & = & 2\langle (c_m'+c_{m-1}\k_m){\bf n}_m-{\bf t}, 
c_1{\bf n}_1+\cdots+c_m{\bf n}_m \rangle \\
\ & = & 2c_m(c_m'+c_{m-1}\k_m); 
\end{array}
$$
$$ 
\begin{array}{rcl}
\mathrm{i.e.}\ \ (R_m^2)' & = & 2c_m(c_m'+c_{m-1}\k_m).
\end{array}
\eqno(3)
$$ 

Thus for $c_m\neq 0$, $c_m'-{(R_m^2)' \over 2c_m}=-c_{m-1}\k_m$. 
This equation together with
the set of equations $(2)$ (using our conventions $c_0=0$ and $c_0'=1$) 
prove Theorem 1. Equation (3) and 
Theorem 3 prove Proposition 0. \fin
\medskip

\noindent
{\bf Proof of Theorem 3 and of its Corollary}.~~By Lemma 2, we have that  
$$C_\g'=(c_m'+c_{m-1}\k_m){\bf n}_m.$$ 
Proposition 1 implies thus that a point of the curve $\g$ is a vertex 
if and only if 
$c_m'+c_{m-1}\k_m=0$. \fin
\medskip

\noindent
{\bf Proof of Theorem 2}.~~The proof will be done by induction. Use
the scalar Frenet equations of Theorem 1 to obtain that
$$\k_1={1 \over c_1},\ \ \k_2={c_1'\over c_2 }={c_1c_1' \over c_1c_2}\ \
\mathrm{and}\ \ 
\k_3={c_2'+c_1\k_2 \over c_3}={c_2'+c_1{c_1'\over c_2 } \over c_3}
={c_2c_2'+c_1c_1' \over c_2c_3}.$$
Suppose that 
$$\k_i={c_{i-1}c_{i-1}'+\cdots+c_2c_2'+c_1c_1' \over c_{i-1}c_i}.
\eqno(4)$$ 
The scalar Frenet equations of Theorem 1 imply that
$c_{i+1}\k_{i+1}=c_i'+c_{i-1}\k_i$. Substitute equation $(4)$ 
to obtain 
$$c_{i+1}\k_{i+1}=c_i'+{c_{i-1}c_{i-1}'+
\cdots+c_2c_2'+c_1c_1' \over c_i}= 
{c_ic_i'+\cdots+c_2c_2'+c_1c_1' \over c_i}.\ \  \square$$

\noindent
{\bf Proof of Theorem 4}.~~We have $R_l^2=c_1^2+\cdots+c_l^2$. Thus 
$R_lR_l'=c_1c_1'+\cdots+c_lc_l'$.  
Combine last equation with the formula of Theorem 2 to obtain
$$R_lR_l'=c_lc_{l+1}\k_{l+1},\ \ \mathrm{for}\ 1\leq l <m.$$
For a generic curve in $\R^{m+1}$ the first $m-1$
Euclidean curvatures are nowhere vanishing and the $m^{th}$ Euclidean curvature may vanish
at isolated points, which do not coincide with the points at which
$R_{m-1}$ is critical. Thus for a generic curve in $\R^{m+1}$, $m>1$,
$R_l'=0$ if and only if either $c_l=0$ or $c_{l+1}=0$ for $1\leq l<m$. 
Moreover, for a smoothly immersed curve the function $c_1=R_1=1/\k_1$ never
vanishes. This proves Theorem 4. \fin
\medskip

\noindent
{\bf Proof of Theorem 5}.~~Write $\s(s)$ for the value of the arc length parameter of 
$C_\g$ at $C_\g(s)$. We assume that the orientations of the parametrised focal 
curve $C_\g$ given by the arc length parameter $s$ of $\g$ and by the arc length 
parameter $\s$ of $C_\g$ coincide. 
Lemma 2 and Theorem 3 imply that at a non-vertex 
of $\g$, the unit tangent vector of the parametrised focal curve $C_\g$ is 
$${\bf T}={(c_m'+c_{m-1}\k_m) \over |c_m'+c_{m-1}\k_m|}{\bf n}_m=\e{\bf n}_m.
\eqno(5)
$$
Moreover, for any non vertex 
$${ds \over d\s}={1 \over |c_m'+c_{m-1}\k_m|}.$$
In order to obtain that 
$${\bf N}_1=\d_1 {\bf n}_{m-1}
\eqno(6)
$$
$$\mbox{and } K_1={|\k_m| \over |c_m'+c_{m-1}\k_m|},$$
derive equation (5) with respect to $\s$ and apply Frenet equations of $\g$ 
taking into account that the first $m-1$ Euclidean curvatures of a generic curve are 
always positive. In the same way, use equation (6) to obtain 
$${\bf N}_2=\d_2 {\bf n}_{m-2}\ \mbox{ and }\ K_2={\k_{m-1} \over |c_m'+c_{m-1}\k_m|}.$$ 
To finish the proof, apply induction process. \fin

\end{document}